   \documentclass[11pt]{amsart}
   \usepackage{latexsym}
   \usepackage{oldlfont}
   \newtheorem{lemma}{Lemma}[section]
   \newtheorem{theorem}[lemma]{Theorem}
   \newtheorem{remark}[lemma]{Remark}
   
   \newtheorem{coro}[lemma]{Corollary}
   \newtheorem{definition}[lemma]{Definition}
   \newcommand{\eps}{\varepsilon}

\newcommand{\e}{\epsilon}

\newcommand{\om}{\omega}
\newcommand{\Om}{\Omega}
\newcommand{\D}{\Delta}
\newcommand{\p}{\partial}
\renewcommand{\phi}{\varphi}

\newcommand{\R}{{\mathbb R}}
\newcommand{\Z}{{\mathbb Z}}

   \input amssym.def
   \input amssym

\title
[Stochastic Two-Layer  Geophysical Flows  ]
{Probabilistic Dynamics\\of  Two-Layer
Geophysical  Flows}

\author{Igor Chueshov}
\address[Igor Chueshov]
{Institute f\"{u}r Dynamische Systeme, FB3\\
Universit\"{a}t  Bremen\\ D-28334 Bremen,
 Germany.  On leave from Department of Mechanics andMathematics,
Kharkov University, 310077 Kharkov, Ukraine}
\email
[Igor D. Chueshov]{chueshov@math.uni-bremen.de}

  \author{Jinqiao Duan}
   \address[Jinqiao Duan]
   {Department of Applied Mathematics\\
   Illinois Institute of Technology\\
   Chicago, IL 60616, USA.}
   \email
   [Jinqiao Duan]{duan@iit.edu}

   \author{Bj{\"o}rn Schmalfu{\ss} }
   \address[Bj{\"o}rn Schmalfu{\ss}]
   {Department of Applied Sciences\\
   University of Technology and Applied Sciences\\
   Geusaer Stra{ss}e\\
   D--06217 Merseburg\\
   Germany\\}
   \email
   [Bj{\"o}rn Schmalfu{\ss}]{schmalfuss@in.fh-merseburg.de}

\date{March 31, 2001}


\begin{document}

\begin{abstract}

The two-layer quasigeostrophic  flow model is an intermidiate system
between the single-layer 2D barotropic flow model and the continuously
stratified, 3D baroclinic flow model. This model is widely used to investigate
basic mechanisms in geophysical flows, such as baroclinic effects,
the Gulf Stream   and subtropical gyres.
The wind forcing acts only on the top
layer. We consider the  two-layer quasigeostrophic  model under stochastic
wind forcing. We first transformed this system into a coupled system of
random partial differential equations and then show that
the asymptotic  probabilistic dynamics of this system depends
only on the top fluid layer.  Namely, in the probability sense and
asymptotically, the dynamics of  the  two-layer quasigeostrophic fluid system
is determinied by the top fluid layer, or, the bottom fluid
layer is slaved by the top fluid layer.
 This conclusion
is true provided that the Wiener process
and  the fluid parameters satisfy a certain condition.
In particular, this latter condition is satisfied when
the trace of the covariance operator of the Wiener process is controled by
a certain upper bound,  and
  the Ekman constant $r$ is sufficiently large.
Note that the generalized time derivative of the Wiener process
models the fluctuating part of the wind stress forcing on the top
fluid layer, and the Ekman constant $r$ measures the   rate for
vorticity decay due to the friction in the bottom Ekman layer.

\end{abstract}

 \maketitle

\newpage
\section{Introduction}

The continuously stratified, three dimensional (3D) baroclinic
quasigeostrophic flow model
describes large scale geophysical fluid motions in the atmosphere and
oceans. This model is much simpler than the primitive flow model
or rotating Navier-Stokes flow model.
When the fluid density  is approximately constant,  this model
reduces to the barotropic, single-layer,  two dimensional (2D)
quasigeostrophic model.
The two-layer quasigeostrophic  flow model,  in which the fluid consists of
two homogeneous fluid layers of  uniform but distinct densities $\rho_1$
and $\rho_2$,
is an intermidiate system
between the single-layer 2D barotropic flow model and the continuously
stratified, 3D baroclinic flow model.

The two-layer quasigeostrophic flow model  has been
used as a theoretical
and   numerical model to understand basic mechanisms,
such as baroclinic effects \cite{Ped87}, wind-drven circulation
\cite{Berloff, Berloff2}, the Gulf Stream \cite{Huang},
fluid stability \cite{Benilov}
and subtropical gyres \cite{Ped96, Oz},  in large scale geophysical flows.

We consider the two-layer quasigeostrophic flow model
(\cite{Ped87}, p. 423):

\begin{equation}\label{q1}
\begin{split}
q_{1t} + J(\psi_1, q_1 + \beta  y )
        &  = \nu \D^2  \psi_1  + f   +  \dot{W}  \; ,\\
q_{2t} + J(\psi_2, q_2 + \beta  y )
       &  = \nu \D^2  \psi_2  -r \D  \psi_2  \; ,
\end{split}
\end{equation}
where potential vorticities $q_1(x, y,t) $, $q_2(x, y,t)$ for
 the
top layer and the bottom layer are
defined via stream functions $\psi_1(x, y,t)$, $\psi_2(x, y,t)$, respectively,
\begin{equation}\label{pv1}
\begin{split}
 q_1 &=  \D \psi_1 - F_1  \cdot (\psi_1 - \psi_2),
\\
 q_2 &=  \D \psi_2 - F_2  \cdot (\psi_2 - \psi_1).
\end{split}
\end{equation}
Here  $x,\,y $ are  Cartesian coordinates in zonal (east),
meridional (north) directions,  respectively;
$(x,y)\in O:=(0,L)\times (0,L)$, where $L$ is a positive number;
 $F_1, F_2$ are positive constants defined by
 (see also \cite{Sal98}, p.87)
\begin{equation*}
\begin{split}
F_1 = \frac{f_0^2 }{gh_1}\frac{\rho_0}{\rho_2-\rho_1},\\
F_2 = \frac{f_0^2 }{gh_2}\frac{\rho_0}{\rho_2-\rho_1},
\end{split}
\end{equation*}
with $g$ the gravitational acceleration;
$h_1, h_2$ the depth of top and bottom layers,
  $\rho_1, \rho_2$ the densities  ($\rho_2> \rho_1$)
of top and bottom layers, respectively;  and
$L, \rho_0$ the characteristic scales for horizontal
length and density of the flows,
 respectively;
 $f_0 +\beta y$ (with $f_0, \beta$ constants) is the Coriolis
parameter and $\beta$ is the meridional gradient of the
Coriolis  parameter, and $ \nu > 0$ is   viscosity..
Note that  $r = f_0 \frac{\delta_E}{2(h_1+h_2)}$ is the  Ekman
constant ( \cite{Ped96}, p.29) which
 measures the intensity of   friction at the bottom boundary layer
(the so-called Ekman layer) or the rate for
vorticity decay due to the friction in the Ekman layer. Here
$\delta_E= \sqrt{2\nu/f_0}$ is the Ekman layer thickness
(\cite{Ped87}, p.188).
   Moreover,  $J(h, g)
=h_xg_y -h_yg_x$ is the Jacobi operator and $\D=\p_{xx} +\p_{yy}$
is the Laplacian operator. Finally, $f(x,y,t)$ is the mean
(deterministic) wind forcing with average zero: $\int_O f dO =0$\\
An important part of the above equation is the white noise term $\dot{W}$,
describing the fluctuating part of the external forcing in the top fluid layer,
see Hasselmann \cite{Has76} and Arnold \cite{Arn00}.
An example of these short  time scale influences is the weather variability
or wind forcing.
A white noise is given as the generalized time derivative of a Wiener
process $W(t)$ in a function space.\\
We assume periodic boundary conditions for
$\psi$ in $x$ and $y$   with period $L$.
In addition,  we impose that
   \[
   \int_O\psi dO      =     0.
   \]
We also assume an appropriate initial condition
\[
q(x, y, 0) = q_0(x,y).
\]
Stochastically forced QGE has been used to investigate various
phenomena in geophysical flows \cite{Holloway, Muller,
Samelson, Griffa, DelSole-Farrell, Brannan}.
Recently Salmon \cite{Salmon}
 introduced some generalized
two-layer ocean  flow models.\\
It is our aim to study the long time dynamics of the stochastic
differential equation with random coefficients but without white noise which we obtain after
a random coordinate transformation.
The solution of this differential equation generates a random dynamical system.
The structure of this differential equation enables us to prove the dissipativity
of the system. In order to investigate the long time dynamics we will apply
the method of {\em determining functionals} in a version for random dynamical systems which
is based on the convergence in probability, see Chueshov et. al \cite{CheDuaSchm00}.
It follows that the asymptotic long time behavior is determined by the asymptotic behavior
of finitely many functionals.
For instance, for these functionals we can choose the Fourier modes given by the eigenfunctions
of the Laplacian. According to this property we can  further
show that  the  functionals
 defined on the top  fluid layer alone can {\em determine} the asymptotical behavior of the complete
two-layer system, when the   fluid parameters and the Wiener
process satisfy a certain condition, such as    the Ekman constant
$r$ is sufficiently large and the trace of the covariance operator
of the Wiener process is controled by some upper bound.
 Note that the
generalized time derivative of the Wiener process models the
fluctuating part of the wind stress forcing on the top fluid
layer, and the Ekman constant $r$ measures the   rate for
vorticity decay due to the friction in the bottom Ekman layer.


We first recall some basic facts in random dynamical systems
in Section 2.  In Section 3, we establish the well-posedness
of the stochastic two-layer quasigeostrophic model by transforming it
into a coupled system of random partial differential equations.
The main results on asymptotic probabilistic determining functionals
are presented in  Sections 4 and 5. Finally we summarize our conclusions
in Section 6.

\section{Random dynamical system and determining functionals}\label{s2}
It is our goal to study the long time dynamics of (\ref{q1}) which is
influenced by random forces. Appropriate tools to treat  this equation
are given by the theory of {\em random dynamical systems}.\\
A random dynamical system   consists
of two components. The first component is a {\em metric dynamical
system} $(\Omega,{\cal F},\Bbb{P},\theta)$ as a model for a
noise, where $(\Omega,{\cal F},\Bbb{P})$ is a probability
space and $\theta$ is a $\cal{F}\otimes
{\cal B}(\Bbb{R}),{\cal F}$ measurable flow: we have
\[
\theta_0={\rm id},\qquad
\theta_{t+\tau}=\theta_t\circ\theta_\tau=:\theta_t\theta_\tau
\]
for $t,\,\tau\in\Bbb{R}$. The measure $\Bbb{P}$ is supposed
to be ergodic with respect to $\theta$. The second  component of a
random dynamical system is a
${\cal B}({\Bbb R}^+)
\otimes{\cal F}\otimes{\cal B}(H),{\cal B}(H)$-measurable
mapping $\phi$ satisfying the {\em cocycle} property
\[
\phi(t+\tau,\omega,x)=\phi(t,\theta_\tau\omega,\phi(\tau,\omega,x)),\qquad \phi(0,\omega,x)=x,
\]
where the phase space $H$ is a separable metric space and $x$ is chosen arbitrarily in $H$. We will denote this
random dynamical system by symbol $\phi$.\\
A standard model for such a noise $\theta$ is the two-sided {\em
Brownian motion}: Let $U$ be a separable Hilbert space. We
consider the probability space
\[
(C_0({\Bbb R},U),{\cal B}(C_0({\Bbb R},U)),\Bbb{P})
\]
where $C_0(\Bbb{R},U)$ is the Fr\'echet space of continuous
functions on $\Bbb{R}$ of uniform convergence on compact intervals which are zero at zero and
${\cal B}(C_0({\Bbb R},U))$ is the corresponding Borel
$\sigma$-algebra. Suppose that we have a covariance operator $Q$
on $U$. Then $\Bbb{P}$ denotes the {\em Wiener measure} with
respect to $Q$. Note that $\Bbb{P}$ is ergodic with respect to
the flow
\begin{equation}\label{NR2}
\theta_t\omega=\omega(\cdot+t)-\omega(t),\qquad\mbox{for
}\omega\in C_0({\Bbb R},U)
\end{equation}
which is called the {\em Wiener shift}.

A main source of a random dynamical system is a random
differential equation.
For example, let us  consider the following evolution equation
 in some Hilbert space
\begin{equation}\label{eq-11}
\frac{du}{d\,t}=F(u,\theta_t\omega),\quad u(0)=x,
\end{equation}
over  some   metric dynamical system $(\Omega,{\cal F},{\Bbb P},\theta)$.
If (\ref{eq-11}) is well-posed for every $\om\in\Om$ and solutions
$u(t,\om;x)$ depends measurably on
$(t,\omega,x)$, then the operator
\[
\phi\; :\; (t,\omega,x)\to u(t,\omega; x)
\]
 defines a random dynamical system (cocycle) $\phi$.
For  detailed presentation of random dynamical systems we
refer to the monograph by L. Arnold \cite{Arn98}.
\medskip\par\noindent

Motivated by deterministic dynamical systems we introduce
several useful notions from the theory of
random dynamical systems.

An $\omega$ depending  closed set $B$ contained in the separable Hilbert space $H$
 is called random
if the mapping $\omega\to\sup_{x\in B(\omega)}\|x-y\|_H$ is a random variable for any $y\in H$.\\
A random dynamical system is called a {\em dissipative}
if there exists  a {\it compact } random set $B$ which
is forward invariant:
\[
\phi(t,\omega,B(\omega))\subset B(\theta_t\omega)
, \; t>0,
\]
and which is absorbing: for  any random
variable $x(\omega)\in H$
there exists a $t_{x}(\omega)>0$ such that if $t\ge t_{x}(\omega)$
\[
\phi(t,\omega,x(\omega))\in B(\theta_t\omega).
\]
For the following we need {\em tempered random variable}:
A random variable $x$ is called tempered if
\[
t\to |x(\theta_t\omega)|
\]
is only subexponentially growing:
\[
\limsup_{t\to\pm\infty}\frac{\log^+|x(\theta_t\omega)|}{|t|}=0
\quad\mbox{a.s.}
\]
This condition is not a very strong restriction because the only alternative is
that the above $\limsup$ is $\infty$ which describes the degenerated case of stationarity,
see Arnold \cite{Arn98}, page 164 f.\\

Our main purpose will be to estimate the degree of freedom of the long time
dynamics of a random dynamical systems which stem from the two-layer
flow problem of the
ocean introduced above.
We will apply the theory of determining functionals to estimate the degree of freedom.
Since dynamics of this model is influenced by random forces we will apply the
theory of determining functionals related to the convergence in probability.\\
We now give our basic definition:

\begin{definition}\label{dfp}
Let $V$ be a Banach space which is continuously embedded into $H$.
Assume that there exists $\tau>0$ such that
$\phi(t,\om,x)\in L^2_{loc}(\tau,+\infty; V)$
for almost all $\om\in\Om$ and $x\in H$.
A set ${\cal L}=\{l_j,\;j=1,\cdots, N\}$ of linear continuous
and linearly independent functionals on $V$ is called
asymptotically determining in probability
if
\[
({\Bbb P}) \lim_{t \to \infty}\int_t^{t+1}
\max_j | l_j(\phi(\tau,\omega,x_1(\omega))-\phi(\tau,\omega,x_2(\omega)))|^2  d \tau  =  0
\]
for two initial conditions $x_1(\omega),\,x_2(\omega)\in H$ implies
$$
({\Bbb P}) \lim_{ t \to \infty}  \|\phi(t,\omega,x_1(\omega))-\phi(t,\omega,x_2(\omega))\|_H
=0.
$$
\end{definition}

The theory of determining functionals was started with the papers
\cite{FP} and \cite{La2} devoted to $2D$ Navier-Stokes equations.
Now this theory is well-developed  for the deterministic
systems (see, e.g., \cite{BC,C3,C4,FMTT,FTiti,JT2} and the references
therein). Some results are also available for
 stochastic systems (see, e.g., \cite{BF, C5,CheDuaSchm00,FL}).
One of the main advantages for this theory is
 the
possibility to localize spatial domains or parameters which are responsible
for long-time dynamics.
The existence of a finite number of determining functionals means that the
long-time behaviour of the system is finite-dimensional. Moreover
the values of these functionals on solutions can be intertreted as a
result of some measurament of the system. From applied point of
view the finiteness of  number of determining functionals means
that we need only a finite number of devices to observe completely all
dynamics of the system.

In this paper we rely on  the following result.

\begin{theorem}\label{t1}
Assume that  a random dynamical system $\phi$ has an
 absorbing  forward invariant random set $B$ in $V$ such that
$\sup_{x\in B(\omega)}\|x\|_V^2$
is bounded by a tempered random variable
and $t\to\sup_{x\in B(\theta_t\omega)}\|x\|_V^2$
is locally integrable.
Let ${\cal L} =\{ l_j : j= 1,...,N\}$ be  a set of linear continuous
and linearly independent functionals on $V$.
Suppose there exist  constants $c_{\cal L}>0$ and  $C_{\cal L}>0$
and a measurable function
$l_{\cal L}(x_1,x_2,\om)$ which maps $V\times V\times\Om$ in $\R$ such that
$l_{\cal L}(x_1,x_2,\om)\ge -c_{\cal L}$ uniformly with respect to
$(x_1,x_2,\om)\in V\times V\times\Om$ and
 for
$x_1(\omega),\,x_2(\omega)\in B(\omega)$ we have
\begin{equation}\label{eqb1}
\begin{split}
V(t,\om)- V(s,\om) & \le
C_{\cal L}\cdot \int_s^t{\cal N}_{\cal L}(\tau,\om)d\tau\\
&+\int_s^t
l_{\cal L}(\phi(\tau,\omega,x_1),\phi(\tau,\omega,x_2),\theta_\tau\omega)
\cdot V(\tau,\om)d\tau
\end{split}
\end{equation}
for all  $t\ge s\ge 0$, where
\begin{align*}
V(t,\om)&=
\|\phi(t,\omega,x_1)-\phi(t,\omega,x_2)\|^2_H,
\\
{\cal N}_{\cal L}(t,\om)&=\max_{j=1,\ldots,N}
|l_j(\phi(t,\omega,x_1)-\phi(t,\omega,x_2))|^2.
\end{align*}
Assume that
\begin{equation*}
\frac{1}{t}{\Bbb E}\left\{\sup_{x_1,x_2\in B(\omega)}\int_0^t
l_{\cal L}(\phi(\tau,\omega,x_1),\phi(\tau,\omega,x_2),\theta_\tau\omega)d\tau\right\}<0
\end{equation*}
for some $t>0$.
Then ${\cal L}$ is a set of asymptotically determining functionals
in probability for random dynamical system  $\phi$.
\end{theorem}
\begin{proof}It is easy to find from relation (\ref{eqb1}) that
\[
V(t,\om)\le V(0,\om)e^{\int_0^tl_{\mathcal L}(s,\om)ds}+
C_{\cal_L}\int_0^t
{\cal N}_{\cal L}(\tau,\om)e^{\int_\tau^tl_{\mathcal L}(s,\om)ds}d\tau,
\]
where
$ l(t,\om)=l_{\cal L}(\phi(t,\omega,x_1),\phi(t,\omega,x_2),\theta_t\omega)$.
Therefore we can apply the argument given in the proof of Theorem~2.2
\cite{CheDuaSchm00}.
\end{proof}

In the next section, we return to the two-layer quasigeostrophic flow model.


\section{Mathematical setup and well-posedness of the random
two-layer  flows }\label{s3}

In the following,  $L^2_{per}$,
$H^s_{per}$ for  $s\in\R$ are the standard Sobolev
spaces of $L$-periodic functions with the zero mean value,
i.e. $\int_O \psi dO = 0$.
Let $(\cdot, \cdot)_0$ and $\| \cdot \|_0 $
denote the standard scalar product and norm in $L^2_{per}$, respectively.
Every element $u(x,y)\in L^2_{per}$ can be represented in
the form
\[
u(x,y)=\sum_{j\in\Z^2, j\neq 0} u_j
L^{-1}\exp\left\{i \frac{2\pi}{L}  (j_1x+j_2y)\right\},
\]
where the Fourier coefficients $u_j$ possess the property
$\bar{u}_j= u_{-j}$  (bar denotes the complex conjugation) and
$$
\|u\|_0^2:=\int_O\vert u\vert^2 dO
\equiv\sum_{j\in\Z^2, j\neq 0} \vert u_j\vert^2 <\infty.
$$
Note that $L^{-1}\exp\left\{i \frac{2\pi}{L}  (j_1x+j_2y)\right\}$ is
 the eigenelement
of the $-\Delta$ with the eigenvalue
$\lambda_1\cdot (j_1^2+j_2^2)$, where
$\lambda_1=\left(\frac{2\pi}{L}\right)^2$ is the smallest eigenvalue.
The norm in $H^s_{per}$ is defined by the formula
$$
\Vert u\Vert^2_s:=\int_O\vert (-\Delta)^{^\frac{s}{2}} u\vert^2 dO
\equiv\lambda_1^{s}
\sum_{j\in\Z^2, j\neq 0}(j_1^2+j_2^2)^s \vert u_j\vert^2.
$$
It is clear that
\begin{equation}\label{n_s}
\Vert\nabla u\Vert^2_s:=\Vert\p_x u\Vert^2_s+\Vert\p_y u\Vert^2_s=
\Vert u\Vert^2_{1+s},\quad s\in\R.
\end{equation}
We also denote ${\bf L}^2_{per}=L^2_{per}\times L^2_{per}$
and  ${\bf H}^s_{per}=H^s_{per}\times H^s_{per}$.
\par
We work on the   phase  space ${\bf H}^{-1}_{per}$ with the scalar product
\[
(q, \bar{q})_{*} =  h_1(\nabla \psi_1,\nabla \bar{\psi}_1)_0
        + h_2 (\nabla \psi_2,\nabla \bar{\psi}_2)_0
    + p (\psi_1-\psi_2,\bar{\psi}_1- \bar{\psi}_2)_0,
\]
where $q=(q_1,q_2)$, $\bar{q}=(\bar{q}_1,\bar{q}_2)$ and
$\psi=(\psi_1,\psi_2)\in {\bf H}_{per}^1$, $\bar{\psi}=(\bar{\psi}_1,\bar{\psi}_2)$.
The relation between  $q$ (resp. $\bar{q}$) and $\psi$ (resp. $\bar{\psi}$)
is defined by (\ref{pv1}).
Here we also use the notation
\[
p = \frac{f_0^2 }{g}\frac{\rho_0}{\rho_2-\rho_1}.
\]
Note that $F_1h_1=F_2h_2=p$. The norm induced by this scalar product
\[
 \|q\|_{*}^2 = (q, q)_{*}= h_1\|\nabla \psi_1\|_0^2
        +   h_2\|\nabla \psi_2\|_0^2
    + p\|\psi_1-\psi_2\|_0^2
\]
is equivalent to the usual norm on ${\bf H}^{-1}_{per}$.
Moreover,  we have the estimate
\begin{equation}
h_1 \|\nabla \psi_1\|_0^2+ h_2 \|\nabla \psi_2\|_0^2\leq
\|q\|_{*}^2  \leq
 a_0\left(h_1\|\nabla \psi_1\|_0^2+ h_2 \|\nabla \psi_2\|_0^2 \right),
\label{equivalent}
\end{equation}
where
\begin{equation*}
a_0=1+\frac{2p}{\lambda_1\min\{h_1, h_2\}}
=1+\frac{2}{\lambda_1}\max\{F_1,F_2\}.
\end{equation*}
To treat the nonlinearity we need the following lemma:

\begin{lemma} \label{Jacobi}
The operator Jacobian verifies, for $u,v,w$ in $H^1_{per}$
\begin{align}
&J(u,v) = - J(v,u),\qquad (J(u,v),v)_0  = 0,\\
&(J(u,v),w)_0 = (J(v,w),u)_0.
\end{align}
Moreover the following estimates hold:
\begin{equation}\label{jac-3}
|( J(u,v),\D u)_0|\le
c_0\|\D v\|_0\cdot\|\nabla u\|_0\cdot \|\D u\|_0,\quad u,\,v\in H_{per}^2;
\end{equation}
\begin{equation}\label{jac-1a}
|( J(u,v),w)_0|\le
c_1\|\D u\|_0\cdot \|\D v\|_0\cdot\|w\|_0,\quad
u,\,v\in H_{per}^2,\,w\in L^2_{per};
\end{equation}
\begin{equation}\label{jac-2a}
|(J(u,v),w)_0|\le
c_1\|\nabla u\|_0\cdot \|\D v\|_0\cdot\|\nabla w\|_0,\quad
u,\,w\in H_{per}^1\,, v\in H_{per}^2.
\end{equation}
Here above
 $c_0=2+(\sqrt{2}\cdot\pi)^{-1}$ and $c_1=c_0\lambda_1^{-\frac{1}{2}}$.
\end{lemma}

\begin{proof} We start with (\ref{jac-3}). It is easy to see that
\begin{align*}
 J(u,v)\cdot\D u=& \frac{1}{2}\cdot v_y\left\{\p_x \left(u_x^2-u_y^2\right)
+2\p_y \left(u_x u_y\right)\right\} \\
+&
 \frac{1}{2}\cdot v_x\left\{\p_y \left(u_x^2-u_y^2\right)
-2\p_x \left(u_x u_y\right)\right\}.
\end{align*}
Consequently relying on property (\ref{n_s}) we have
\begin{align*}
|( J(u,v),\D u)_0|\le
& \frac{1}{2}\cdot \|v_y\|_1 \left(\|\p_x \left(u_x^2-u_y^2\right)\|_{-1}
+2\|\p_y \left(u_x u_y\right)\|_{-1}\right) \\
+&
\frac{1}{2}\cdot \|v_x\|_1 \left(\|\p_y \left(u_x^2-u_y^2\right)\|_{-1}
+2\|\p_x \left(u_x u_y\right)\|_{-1}\right) \\
\le &
 \frac{1}{2}\cdot\left( \|v_x\|_1 + \|v_y\|_1 \right)\cdot
\left(\| u_x\|_{L^4}^2+\| u_y\|_{L^4}^2
+2\|u_x u_y\|_0\right)\\
\le
&
\left( \|v_x\|_1 + \|v_y\|_1 \right)\cdot
\left(\| u_x\|_{L^4}^2+\| u_y\|_{L^4}^2 \right) \\
\le
&
\sqrt{2}\cdot \|\D v\|_0
\cdot\left(\| u_x\|_{L^4}^2+\| u_y\|_{L^4}^2 \right).
\end{align*}
Using the inequality (see, e.g., \cite{FMTT})
\begin{equation*}
\Vert u \Vert_{L^4}\le a_1\cdot \Vert u \Vert^{1/2}
\cdot \Vert \nabla u \Vert^{1/2}\quad\mbox{with}\quad
a_1=\left( ( 2\pi)^{-1}+\sqrt{2}\right)^{1/2},
\end{equation*}
we obtain  (\ref{jac-3}).

In a similar way we can establish the estimates
\begin{equation*}
|(J(u,v),w)_0|\le
c_0\|\nabla u\|_0^{1/2}\cdot \|\D u\|_0^{1/2}\cdot
\|\nabla v\|_0^{1/2}\cdot \|\D v\|_0^{1/2}\cdot\|w\|_0,
\end{equation*}
\begin{equation*}
|(J(u,v),w)_0|\le
c_0\| u\|_0^{1/2}\cdot \|\nabla u\|_0^{1/2}\cdot
\|\nabla v\|_0^{1/2}\cdot \|\D v\|_0^{1/2}\cdot\|\nabla w\|_0,
\end{equation*}
which easily imply  (\ref{jac-1a}) and (\ref{jac-2a}).
\end{proof}

We now transform the two-layer model  (\ref{q1})
containing white noise into a system of coupled
random partial differential equations.
The reason for  taking
such a transformation is that we need some particular a priori estimates
for the trajectories of the system. Often these a priori estimates can be calculated
by the Gronwall lemma.
For equations containing a white noise there exists no Gronwall
lemma. But we can use this technique for the transformed
random partial differential equations.\\
For this purpose we
introduce an Ornstein-Uhlenbeck process $\eta (x,y,t, \omega)$
in $L^2_{per}$. This process is defined  by
the solution of the following linear stochastic partial differential equation
\begin{equation}\label{OU-pr}
\eta_{t}  =  \nu (k+1)\D \eta + \dot{W},
\end{equation}
with periodic boundary condition and with some initial condition $x\in L_{per}^2$
where $k>0$ is a free {\em control} parameter
and $W$ is a Wiener process
in $L^2_{per}$. We suppose that the covariance operator $Q$ of this
Wiener process has a finite trace. As mentioned in Section \ref{s2} such a Wiener process
generates a metric dynamical system denoted by $(\Omega,\mathcal{F},\mathbb{P},\theta)$
where $\theta$ is the Wiener shift introduced in (\ref{NR2}).
It is well known that
this equation has a stationary solution which is generated by
a {\em Gaussian} random
variable $\eta$ in $H_{per}^1$. In particular,
the mapping
\[
t\to\eta(\theta_t\omega)\in L_{loc}^2(-\infty,\infty;H_{per}^1)
\]
solves this equation. For  moments of $\eta$ we obtain:
\begin{equation}\label{constanten}
{\mathbb E}\|\eta\|_{1}^2\le \frac{{\rm tr_0}Q}{2\nu(k+1)},\quad
{\mathbb E}\|\eta\|_{1}^{2n} \le C_n\left(\frac{{\rm
tr_0}Q}{\nu(k+1)}\right)^n,\; n \in{\mathbb N},\;C_n>0.
\end{equation}
We introduce new variables
\begin{equation}\label{zz}
\tilde{q}_1:=q_1-\eta,\; q_2,\;  \tilde{\psi}_1:=\psi_1+\xi_1,\;
 \tilde{\psi}_2:=\psi_2+\xi_2,
\end{equation}
where stationary process $\eta$ solves the problem (\ref{OU-pr})
and $\xi_1$ and $\xi_2$  are defined such that the elliptic
equations (\ref{pv1}) remain  of the same form
\begin{equation}\label{pv1-new}
\begin{split}
\tilde{q}_1 &= \D \tilde{\psi}_1 - F_1 \cdot (\tilde{\psi}_1 - \tilde{\psi}_2),
\\
 q_2 &=  \D \tilde{\psi}_2 - F_2  \cdot (\tilde{\psi}_2 - \tilde{\psi}_1),
\end{split}
\end{equation}
The processes $\xi_1$ and  $\xi_2$ are solutions of  the linear
elliptic equations
\begin{equation}\label{pv1-new1}
\begin{split}
 \D \xi_1 - F_1 \cdot (\xi_1 - \xi_2) =& -\eta,
\\
 \D \xi_2 - F_2 \cdot (\xi_2 - \xi_1)   =&    0
\end{split}
\end{equation}
and can be presented in the form
\begin{align*}
\xi_1&=\frac{1}{F_1+F_2}\left\{
F_2 (-\D)^{-1}+F_1(-\D+ F_1+F_2)^{-1}\right\}\eta
\\
\xi_2&=\frac{F_2}{F_1+F_2} \left\{
(-\D)^{-1}-(-\D+  F_1+F_2)^{-1}\right\}\eta.
\end{align*}
Thus the  processes $\xi_1$ and $\xi_2$ are   smoother in spatial variables
than $\eta$.  In fact after simple calculations we have the estimates
\begin{equation}
\|\xi_i\|_{s+2}      \leq    \|\eta\|_s,
\quad i=1,2,\;\quad\text{and}\quad
\|\xi_1-\xi_2\|_{s+2}\le \|\eta\|_s\quad  s\in\R.
\label{xi}
\end{equation}
Now we obtain the  coupled random partial differential equations
for  new potential vorticities $\tilde{q}_1$, $q_2$:
\begin{align*}
\tilde{q}_{1t} + &J(\tilde{\psi}_1-\xi_1, \tilde{q}_1+\eta + \beta  y )
          = \nu \D^2  \tilde{\psi}_1  + f  \nonumber \\
&-
\nu\D^2\xi_1- \nu (k+1)\D \eta  \; ,
\\
q_{2t} + &J(\tilde{\psi}_2-\xi_2, q_2 + \beta  y )
         = \nu \D^2  \tilde{\psi}_2  -r \D  \tilde{\psi}_2  \nonumber\\
&-\nu \D^2  \xi_2  + r \D  \xi_2 .
\end{align*}
We treat $\eta, \xi_1, \xi_2$
as known processes. Using (\ref{pv1-new1}) we have
\begin{align*}
-\nu\Delta^2\xi_1-\nu(k+1) \Delta\eta
      =&-\nu F_1(\Delta\xi_1-\Delta\xi_2)-\nu k\Delta\eta,\\
-\nu\Delta^2\xi_2=&-\nu F_2(\Delta\xi_2-\Delta\xi_1).
\end{align*}
For convenience,  we drop the tilde  and rewrite the above system.
Thus we finally get the coupled system of
random  partial differential equations
\begin{equation}\label{Eqn1}
\begin{split}
q_{1t} +J(\psi_1-\xi_1, q_1+\eta + \beta  y ) = &  \nu \D^2 \psi_1  + f\\
&-\nu F_1(\Delta\xi_1-\Delta\xi_2)-\nu k\Delta\eta\\
q_{2t} +  J(\psi_2-\xi_2, q_2 + \beta  y ) = & \nu \D^2 \psi_2  -r \D  \psi_2
\\
& -\nu F_2(\Delta\xi_2-\Delta\xi_1)   + r \D  \xi_2\\
\end{split}
\end{equation}
with
\begin{equation}\label{Eqn3}
\begin{split}
 q_1 =&  \D \psi_1 - F_1  \cdot (\psi_1 - \psi_2),
\\
 q_2 =&  \D \psi_2 - F_2  \cdot (\psi_2 - \psi_1),
\end{split}
\end{equation}
in the class of $L$-periodic functions with initial data
$q(x,y,0)=q_0(x,y)\equiv (q_{01}(x,y),q_{02}(x,y))\in {\bf H}_{per}^{-1}$,
where $\eta$ is the stationary solution to (\ref{OU-pr})
and $\xi_1$ and $\xi_2$ are solves (\ref{pv1-new1}) in $H^2_{per}$.

For the rest of the paper, we work on this coupled system of
random partial differential
equations for stochastically forced two-layer quasigeostrophic
fluid system.

The coefficients of  coupled system (\ref{Eqn1}), (\ref{Eqn3})
have similar properties
as the coefficients of the corresponding deterministic
two-layer quasigeostrophic  system (cf. \cite{Bernier94}, where $N$-layer
model with other boundary conditions is considered). Therefore,
similar to  \cite{Bernier94}, using the Galerkin method and the compactness
argument we can prove the following assertion on
 well-posedness of problem (\ref{Eqn1}), (\ref{Eqn3})  in the class
of $L$-periodic functions.

\begin{theorem}\label{t0}
Let $q_{0} \in {\bf H}^{-1}_{per}$ and $f \in L^2_{per}$. Then for all
$\om\in\Om$ and for all $T > 0$,
the system (\ref{Eqn1}), (\ref{Eqn3})  has a
unique solution
$\{q(t), \psi(t)\}$ such that
\[
q\in C([0,T];{\bf H}_{per}^{-1})\cap L^2(0,T;{\bf L}_{per}^2)\cap
L^{2}_{loc}(0,T;{\bf H}^{1}_{per}).
\]
The function $\psi$ associated to
$q$ by (\ref{Eqn3}) satisfies
\[
\psi \in
C([0,T];{\bf H}^{1}_{per})\cap L^2(0,T;{\bf H}^2_{per})\cap
L^{2}_{loc}(0,T;{\bf H}^{3}_{per}).
\]
The solution depends continuously on the initial condition
$q_0\in {\bf H}_{per}^{-1}$.
\end{theorem}
By the uniqueness assertion of the last Theorem the solution
$t\to q(t)$ generates a random dynamical system
$\phi$ with continuous ${\bf H}_{per}^{-1}\ni q\to\phi(t,\omega,q)$ on the
phase space ${\bf H}_{per}^{-1}$.

\section{Dissipativity of the random dynamical system}

         Dynamical systems generated by many nonlinear parabolic
differential equations have the  dissipative property  which
means that there exists a compact set  absorbing the states of
the system in finite time. Usually one can choose such a set which is also
forward invariant.\\
We now construct  an absorbing forward invariant set for the
random dynamical system  generated by (\ref{Eqn1}). This set will be  a random set.
\begin{theorem}
There exists a compact random set
$B(\omega)\subset {\bf H}_{per}^{-1}$ such that
\begin{equation}\label{eqn3}
\begin{split}
&\phi(t,\omega,B(\omega))\subset B(\theta_t\omega)\quad \text{for }t\ge 0,\\
&\phi(t,\omega,q(\omega))\subset B(\theta_t\omega)
\quad\text{for }t\ge t_0(\omega,q)
\end{split}
\end{equation}
where $q$ is a  random variable with values in ${\bf H}_{per}^{-1}$.
\end{theorem}
We now divide the proof of this theorem in some lemmata. We start
with the following:
\begin{lemma}\label{l0}
Let $q(t)$ be the solution of (\ref{Eqn1}). Then $q(t)$ satisfies the following inequality
\begin{align*}
\frac{d}{dt}\|q(t)\|_{\ast}^2+&
\nu (h_1 \|\Delta\psi_1(t)\|_0^2+h_2 \|\Delta\psi_2(t)\|_0^2)
 \\
\le& d_0\cdot \|\eta (\theta_t\omega)\|^2_0 \cdot
(h_1 \|\nabla\psi_1\|_0^2+h_2 \|\nabla\psi_2\|_0^2)
+m(\theta_t\omega),
\end{align*}
where
\[
m(\omega)= d_1\|\eta(\omega)\|_0^4+d_2\|\eta(\omega)\|_0^2 +d_3
\]
and
\begin{align*}
d_0=&\frac{6c_0^2}{\nu}\left( 1+ \frac{ p^2\nu }{\lambda_1^2
\min\{h_1,h_2\}}\right), \\
d_1=&\frac{6c_0^2h_1}{\nu\lambda_1}, \\
d_2=& 9\left(
\frac{\beta^2(h_1+h_2)}{\nu\lambda_1^3}+
\frac{\nu p^2}{\lambda_1^2}\left(\frac{1}{h_1}+\frac{1}{5h_2}\right)+
\frac{rh_2}{18\lambda_1} +k^2\nu h_1
\right)\;,
\\
d_3=& \frac{9h_1}{\nu\lambda_1}\|f\|_{-1}^2.
\end{align*}
\end{lemma}
\begin{proof}
Multiplying  the equations (\ref{Eqn1})  by
$-2h_1\psi_1,\,-2h_2\psi_2$, respectively,
 and then adding together,  we obtain
\begin{align*}
\frac{d}{dt}\|q(t)\|_{\ast}^2+&
2\nu (h_1 \|\Delta\psi_1\|_0^2+h_2 \|\Delta\psi_2\|_0^2)
+2rh_2 \|\nabla \psi_2\|_0^2 \\
=& 2h_1 (J(\psi_1-\xi_1,q_1+\eta+\beta y),\psi_1)_0 \\
&+2h_2 (J(\psi_2-\xi_2,q_2+\beta y),\psi_2)_0
-2h_1  (f,\psi_1)_0+2h_1 \nu k (\Delta\eta, \psi_1)_0 \\
&+2 \nu p(\Delta\xi_1-\Delta\xi_2, \psi_1 )_0
 + 2 \nu p(\Delta\xi_2-\Delta\xi_1,   \psi_2 )_0 \\
&-2 r h_2 (\Delta\xi_2, \psi_2)_0.
\end{align*}
We have by Lemma \ref{Jacobi}
\[
(J(\psi_1,q_1+\eta),\psi_1)_0=0, \; (J(\psi_2, q_2),\psi_2)_0=0.
\]
A simple calculation
shows   that $(J(\psi_i,\beta y),\psi_i)_0=0$.
We  now estimate the right hand side terms.
On account of (\ref{Eqn3}), (\ref{jac-3}), (\ref{jac-2a}), the bilinearity of $J$ and
(\ref{xi}) we can estimate
\begin{align*}
-2\sum_{i=1}^2&h_i(J(\xi_i,q_i),\psi_i)_0=-2\sum_{i=1}^2h_i(J(\xi_i,\Delta \psi_i),\psi_i)_0
+2p(J(\xi_1-\xi_2,\psi_2),\psi_1)_0\\
\le &2c_0\sum_{i=1}^2h_i\|\Delta\xi_i\|_0\|\Delta\psi_i\|_0\|\nabla\psi_i\|_0
+2c_1 p\|\Delta(\xi_1-\xi_2)\|_0\|\nabla\psi_1\|_0\|\nabla\psi_2\|_0\\
\le&\frac{\nu}{6}(h_1\|\Delta\psi_1\|_0^2+h_2\|\Delta\psi_2\|_0^2)+
6\frac{c_0^2}{\nu}\|\eta\|_0^2(h_1\|\nabla\psi_1\|_0^2+h_2\|\nabla\psi_2\|_0^2)
\\
&+
\frac{c_1 p\|\eta\|_0}{\min\{h_1,h_2\}}(h_1\|\nabla\psi_1\|_0^2+h_2\|\nabla\psi_2\|_0^2) \\
\le&\frac{\nu}{3}(h_1\|\Delta\psi_1\|_0^2+h_2\|\Delta\psi_2\|_0^2)+
d_0\|\eta\|_0^2(h_1\|\nabla\psi_1\|_0^2+h_2\|\nabla\psi_2\|_0^2).
\end{align*}
Similarly, due to (\ref{jac-1a}) and (\ref{xi})
\[
-2h_1(J(\xi_1,\eta),\psi_1)_0\le
2c_1 h_1\|\Delta\xi_1\|_0\|\eta\|_0\|\D\psi_1\|_0
\le \frac{6c_1^2h_1\|\eta\|_0^4}{\nu}+\frac{\nu h_1}{6}\|\D\psi_1\|^2_0.
\]
By the Cauchy-Schwarz inequality
\begin{align*}
&-2\sum_{i=1}^2h_i(J(\xi_i,\beta y),\psi_i)_0\le
2\beta \sum_{i=1}^2 h_i\|\nabla\xi_i\|_{-1}\|\nabla\psi_i\|_0\\
&
\le
\frac{9\beta^2(h_1+h_2)}{\nu\lambda_1^3}\|\eta\|_0^2
+\frac{\nu}{9}(h_1\|\D\psi_1\|_0^2+h_2\|\D\psi_2\|_0^2),\\
&
-2h_1(f,\psi_1)_0\le 2h_1\|f\|_{-1}\|\nabla\psi_1\|_0\le
\frac{9 h_1}{\nu\lambda_1}\|f\|_{-1}^2+\frac{\nu h_1}{9}\|\D\psi_1\|_0^2\\
&
+2h_1\nu(k\Delta \eta,\psi_1)_0\le 9k^2\nu h_1\|\eta\|_0^2
+\frac{\nu h_1}{9}\|\D\psi_1\|_0^2.
\end{align*}
Using (\ref{xi}) once more
\begin{align*}
2\nu p(\Delta(\xi_1-\xi_2),\psi_1)_0&\le
2\nu p\|\Delta\xi_1-\Delta\xi_2\|_{-2}\|\D\psi_1\|_0\\
&\le
\frac{9\nu p^2}{\lambda_1^2h_1}\|\eta\|_0^2+
\frac{\nu h_1}{9}\|\D\psi_1\|_0^2,\\
-2\nu p(\Delta(\xi_1-\xi_2),\psi_2)_0
&\le
\frac{9\nu p^2}{5\lambda_1^2h_2}\|\eta\|_0^2+\frac{5\nu h_2}{9}
\|\D\psi_2\|_0^2,\\
-2rh_2(\Delta\xi_2,\psi_2)_0&
\le
\frac{rh_2}{2 \lambda_1}\|\eta\|_0^2+2rh_2\|\nabla\psi_2\|_0^2.
\end{align*}
Adding all these inequalities together we obtain the conclusion.
\end{proof}

We now consider the random variable $\eta$ defined in Section \ref{s2}.
Recall that $\eta$
depends on the control parameter $k$. If $k$ is chosen large enough then
particular moments of $\eta$ are small. Especially we can formulate:
\begin{lemma}\label{l1}
Let $W$ be a Wiener process in $L^2_{per}$ with finite trace of the covariance.
Then under assumptions
\begin{equation}\label{assump}
\frac{2d_0 a_0 {\rm tr}_0Q}{\lambda_1^2\nu^2 (k+1)}<1,\quad
\frac{16d_0  {\rm tr}_0Q}{\lambda_1^2\nu^2 (k+1)^2}<1
\end{equation}
the random variable
\[
R_0(\omega):=\int_{-\infty}^0
e^{\frac{\nu\lambda\tau}{a_0}\tau+d_0\int_{\tau}^0
\|\eta(\theta_{\tau^\prime}\omega)\|_0^2d\tau^\prime}
m(\theta_\tau\omega)d\tau
\]
 is finite and tempered. Moreover
 \begin{equation*}
\left({\mathbb E} R_0^2\right)^{1/2}\le
d_4^2 \left(\frac{3a_0}{2\nu\lambda_1}\right)^{3/2}\cdot
\left(
\frac{2\lambda_1\nu}{a_0}-
\frac{4d_0  {\rm tr}_0Q}{\lambda_1\nu (k+1)}
\right)^{-1/2},
\end{equation*}
where
\[
d_4=  C_8^{1/4} \frac{d_1  ({\rm tr}_0Q)^2}{\lambda_1^2\nu^2 (k+1)^2}
+C_4^{1/4} \frac{d_2  {\rm tr}_0Q}{\lambda_1\nu (k+1)}+d_3
\]
is an estimate for $({\mathbb{E}}m^4)^{1/4}$ (the constants
$C_8,\,C_4$ are defined in (\ref{constanten})).
\end{lemma}
The proof of this lemma can be found in Chueshov et. al.
\cite{CheDuaSchm00} for an Ornstein-Uhlenbeck process in another Hilbert space.
However the argument given there is of a general nature.
\\
We now construct a set satisfying (\ref{eqn3}).
\begin{lemma}\label{l5}
Let $R(\omega):=a R_0(\omega)$ for some $a>1$ and $R_0$ as in Lemma \ref{l1}.
Then the closed  ${\bf H}_{per}^{-1}$-ball $B(0,R(\omega)^\frac{1}{2})$ fullfills (\ref{eqn3}) provided conditions (\ref{assump}) hold.
\end{lemma}
\begin{proof}
Using Lemma~\ref{l0} and relation (\ref{equivalent}) we have
\[
\frac{d}{dt}\|q(t)\|_{\ast}^2\le
\left(-\frac{\nu\lambda_1}{a_0}+  d_0\cdot \|\eta (\theta_t\omega)\|^2_0
\right)\cdot\|q(t)\|_{\ast}^2+m(\theta_t\omega).
\]
Let $q_0=q(0)$ and $\rho(t,\omega, \|q_0\|_\ast^2)$
be the solution of one dimensional
random affine equation
\begin{equation}\label{eq4}
\frac{d\rho(t)}{dt}+\frac{\nu\lambda_1}{a_0}\rho=
d_0\|\eta(\theta_t\omega)\|_0^2\rho+m(\theta_t\omega),
\quad
\rho(0,\omega,\|q_0\|_\ast^2)=\|q_0\|_\ast^2.
\end{equation}
A comparison argument gives that
\[
\|\phi(t,\omega,q_0)\|_\ast^2\equiv\| q(t)\|_\ast^2
\le \rho(t,\omega,\|q_0\|_\ast^2).
\]
Here $\phi$ is the dynamical system introduced in Section \ref{s3}:
$\phi(t,\omega,q_0)=q(t)$, where $q(t)$ is the solution to (\ref{Eqn1}) with
the initial data $q_0$.
Equation
(\ref{eq4}) has the stationary solution given by $t\to R_0(\theta_t\omega)$:
\[
\rho(t,\omega,R_0(\omega))=R_0(\theta_t\omega).
\]
This can be checked by the variation of constants formula.
This solution is exponentially
attracting which follows from the variation of constants formula again:
\begin{align*}
|R_0(\theta_t\omega)-\rho(t,\omega,\|q_0\|_\ast^2)|
&=|\rho(t,\omega,R_0(\omega))-\rho(t,\omega,\|q_0\|_\ast^2)|\\
&\le
e^{\int_0^t\left(d_0\|\eta(\theta_\tau\omega)\|^2_0-
\frac{\nu\lambda_1}{a_0}\right)d\tau}
\left( R(\omega)+\|q_0\|_\ast^2\right)
\end{align*}
which tends to zero exponentially fast. Indeed,
it follows from (\ref{constanten}) and (\ref{assump}) that
for a sufficient small $\eps>0$,
\[
\int_0^td_0\|\eta(\theta_\tau\omega)\|_0^2d\tau<\frac{\nu\lambda_1-\eps}{a_0}t,\qquad
\frac{\nu\lambda_1-\eps}{a_0}>0
\]
for large $t>0$ under conditions (\ref{assump}).
\end{proof}
It remains to prove the existence  of a compact set  $B$ satisfying
(\ref{eq4}).
\begin{lemma}
Suppose that the random variable $R(\omega)$ is defined in Lemma \ref{l5}.
The set
\[
B(\omega):=\overline{\phi(1,\theta_{-1}\omega,B(0,R(\theta_{-1}\omega)^\frac{1}{2}))}
\]
is a compact absorbing forward invariant random set. Moreover
\begin{equation}\label{tempered}
\om\mapsto\sup\left\{ \|\D\psi_1\|_0^2+\|\D\psi_1\|^2_0\; :\;
(q_1, q_2)\in B(\om)\right\}
\end{equation}
is a tempered random variable ($\psi_1$ and $\psi_2$ are defined by
(\ref{Eqn3})).
\end{lemma}
\begin{proof}
Since $\phi(t,\omega,\cdot)$ is completely continuous for $t>0$ (see
the regularity assertion of Theorem \ref{t0}) the sets
$B(\omega)$ are compact. Since $R$ is a random variable the ball $B(0,R^\frac{1}{2})$ is a
random set. The continuity  of $\phi(t,\omega,\cdot)$ allows us to
conclude that $B$ is a random set. The construction of $B$ ensures
that that set in absorbing and forward invariant.
The temperedness of (\ref{tempered}) can be proved in the same way as in
\cite{CheDuaSchm00} for the $2D$ Navier-Stokes equations.
\end{proof}

Hence we have shown that the two-layer quasigeostrophic flow
system is dissipative in the probabilistic sense.
In the next section we consider asymptotic probabilistic
determining functionals.

\section{Determinging functionals}

 In this section,   we compare two solutions
as a way to understand asymptotic dynamics in the probabilistic sense.
Consider the difference of two solutions  $\hat{q}$ and $\bar{q}$,
corresponding to  the stream functions
$ \hat{\psi}$,  $ \bar{\psi}$ .
We set
$$ q=(q_1,q_2)=\hat{q} - \bar{q},$$
corresponding to the stream function
$$\psi =\hat{\psi} - \bar{\psi}.$$


We get random partial differential
equations for the solution difference  $q=\hat{q} - \bar{q}$ from
(\ref{Eqn1}):
 \begin{align}
q_{1t} = &  \nu \D^2 \psi_1
            \nonumber \\
&-J(\psi_1, \hat{q}_1 + \beta  y ) -J(\bar{\psi}_1, q_1)
 -J(\psi_1, \eta) + J(\xi_1, q_1)  \; ,
\label{QQ}\\
q_{2t}  = & \nu \D^2 \psi_2  -r \D  \psi_2
                     \nonumber\\
&-  J(\psi_2, \hat{q}_2 + \beta  y )-J(\bar{\psi}_2, q_2)
 + J(\xi_2, q_2) \; ,
\label{QQQ}
\end{align}

Multiplying (\ref{QQ}) by $-h_1 \psi_1$, and (\ref{QQQ}) by $-h_2 \psi_2$
and adding together, we get
\begin{equation}\label{estimate1}
\begin{split}
 \frac12 \frac{d}{dt}\|q\|_{*}^2  = & -\nu [h_1 \|\D \psi_1\|_0^2
+h_2 \|\D \psi_2\|_0^2] - r h_2 \| \nabla \psi_2\|_0^2  \\
+&  h_1 (J(\bar{\psi}_1-\xi_1, q_1), \psi_1)_0
+ h_2 (J(\bar{\psi}_2-\xi_2, q_2), \psi_2)_0.
\end{split}
\end{equation}
We start with
\[
h_1(J(\bar{\psi}_1-\xi_1, q_1), \psi_1)_0\equiv I_1+I_2,
\]
where
\[
I_1= h_1(J(\bar{\psi}_1-\xi_1, \D\psi_1), \psi_1)_0,\quad
I_2= p(J(\bar{\psi}_1-\xi_1, \psi_2), \psi_1)_0.
\]
By (\ref{jac-3}) we have
\[
| I_1|\le c_0 h_1  \|\bar{\psi}_1-\xi_1\|_2 \|\D\psi_1\|_0 \|\nabla\psi_1\|_0.
\]
Using (\ref{jac-2a}) we obtain
\[
| I_2|\le c_1p
\|\bar{\psi}_1-\xi_1\|_2 \|\nabla\psi_2\|_0 \|\nabla\psi_1\|_0.
\]
Thus
\[
|I_2|\le c_1p\lambda_1^{-\frac{1}{2}}
\|\bar{\psi}_1-\xi_1\|_2 \|\D\psi_2\|_0 \|\nabla\psi_1\|_0.
\]
Using the inequality
$ab  \leq  \e a^2 +\frac1{4\e}b^2$,  we get
\begin{equation}\label{int-h-1}
\begin{split}
h_1(J(\bar{\psi}_1-\xi_1,& q_1), \psi_1)_0 \le
\frac{\nu}{4}\left( h_1\|\D\psi_1\|_0^2+ h_2\|\D\psi_2\|_0^2\right)\\
& +
\frac{h_1c_0^2}{\nu}(1+ F_1 F_2 \lambda_1^{-2})
\|\bar{\psi}_1-\xi_1\|_2^2  \|\nabla\psi_1\|_0^2.
\end{split}
\end{equation}
In a similar way we have
\begin{equation}\label{int-h-2}
\begin{split}
h_2(J(\bar{\psi}_2-\xi_2&, q_2), \psi_2)_0 \le
\frac{\nu}{4}\left( h_1\|\D\psi_1\|_0^2+ h_2\|\D\psi_2\|_0^2\right)
\\
& +
\frac{h_2c_0^2}{\nu}(1+ F_1 F_2 \lambda_1^{-2})
\|\bar{\psi}_2-\xi_2\|_2^2  \|\nabla\psi_2\|_0^2.
\end{split}
\end{equation}
Using (\ref{estimate1}), (\ref{int-h-1}) and
(\ref{int-h-2}) we obtain
\begin{equation}\label{estimate1a}
\begin{split}
 \frac{d}{dt}\|q\|_{*}^2 \le & -\nu [h_1 \|\D \psi_1\|_0^2
+h_2 \|\D \psi_2\|_0^2] - 2r h_2 \| \nabla \psi_2\|_0^2   \\
&+
b_0
\left( h_1
\|\bar{\psi}_1-\xi_1\|_2^2  \|\nabla\psi_1\|_0^2
+
 h_2
\|\bar{\psi}_2-\xi_2\|_2^2  \|\nabla\psi_2\|_0^2\right)\;,
\end{split}
\end{equation}
where
$b_0=\frac{2c_0^2}{\nu}
(1+ F_1 F_2\lambda_1^{-2})$.
This estimate is the main point in the construction of determining
functionals.
It is our aim to obtain an integral equation like (\ref{eqb1}).
To this end
we first consider the case when parameter $r>0$ is arbitrary.
Let
 ${\cal L}$ = $\left\{l_j \right\}_1^N $ be a set of linearly independent
bounded linear functionals on the space ${\bf H}^2_{per}$.
 Assume that the set  ${\cal L}$ possesses the property
\begin{equation}
\label{2dieze}
\|\psi\|_{1} \leq  C_{\cal L} \max_j |l_j(\psi)| +
\varepsilon_{\cal L} \|\psi\|_{2},
\end{equation}
for any $\psi \in {\bf H}^2_{per}$ with some positive
constant $ C_{\cal L}$ and $\varepsilon_{\cal L}$.
We note that the best possible value of the parameter $\varepsilon_{\cal L}$
is called the completeness defect of the family ${\cal L}$ with respect
of the pair of the spaces ${\bf H}^1_{per}$ and ${\bf H}^2_{per}$
(see \cite{C3,C4}) and the smallness of the parameter $\varepsilon_{\cal L}$
is crucial in the subsequent considerations. We refer to
\cite{C3,C4} for general properties of completeness defect and for
estimates of $\varepsilon_{\cal L}$ for several families of functionals
on Sobolev spaces.

From (\ref{2dieze})  for
$\psi=(\sqrt{h_1}\psi_1,\sqrt{h_2}\psi_2)$ we have
\begin{align*}
\left(h_1\|\D\psi_1\|_0^2+h_2\|\D\psi_2\|_0^2\right) \ge&
\frac{1-\delta}{\eps_{\cal L}^2}\left(h_1\|\nabla \psi_1\|_0^2+
\|h_2\nabla \psi_2\|_0^2\right)\\
&-C_{\delta,{\cal L}}\max_{j=1,\cdots,N}
|l_j(\sqrt{h_1}\psi_1,\sqrt{h_2}\psi_2)|^2.
\end{align*}
for any $0<\delta<1$ with appropriate positive constant
$C_{\delta,{\cal L}}$. Therefore from (\ref{estimate1a}) and
(\ref{equivalent}) we obtain
\begin{equation}
\label{el-1}
 \frac{d}{dt}\|q\|_{*}^2 \le l^{(1)}_{\cal L}(\bar{\psi}, \om)
\|q\|_{*}^2 + C_{\delta,{\cal L}}\max_{j=1,\cdots,N}
|l_j(\sqrt{h_1}\psi_1,\sqrt{h_2}\psi_2)|^2,
\end{equation}
where
\begin{equation*}
 l^{(1)}_{\cal L}(\bar{\psi}, \om)=
-\frac{1-\delta}{a_0\eps_{\cal L}^2}+
b_0\max\left\{\|\bar{\psi}_1-\xi_1\|_2^2\;,
\|\bar{\psi}_2-\xi_2\|_2^2\right\}.
\end{equation*}
The application of Theorem \ref{t1} gives the following assertion.

\begin{theorem}  \label{main-1}
Let    ${\cal L} =  \{ l_j : j= 1,\cdots,N\}$ be a finite set of
linearly independent continuous functionals on the space ${\bf H}^2_{per}$.
We assume that this set satisfies (\ref{2dieze}). Let (\ref{assump}) hold.
If
\begin{equation*}
\eps_{\cal L}<
\frac{\nu\sqrt{\min\{h_1,h_2\}}}{\sqrt{2a_0b_0\Sigma}},
\end{equation*}
where
\begin{equation}\label{si-2}
\Sigma=d_0({\mathbb E}\|\eta\|_0^4)^\frac{1}{2}({\mathbb E}R^2)^\frac12+
{\mathbb E}m
+\nu\min\{h_1,h_2\}{\mathbb E}\|\eta\|_0^2
\end{equation}
and $m$, involving the Wiener process through
the  Ornstein-Uhlenbeck process $\eta$, is defined in Lemma \ref{l0},
then ${\cal L}^{h_1,h_2}=\{l_j^{h_1,h_2}:j=1,\cdots, N\}$
where $l_j^{h_1,h_2}(\psi_1,\psi_2):=l_j(\sqrt{h_1}\psi_1,\sqrt{h_2}\psi_2),\;l_j\in {\cal L}$
is a set of asymptotically determining functionals in
probability for the  stochastically forced two-layer
quasigeostrophic fluid system (\ref{Eqn1}), (\ref{Eqn3}).
\end{theorem}
\begin{proof}
We integrate (\ref{el-1}).
The conditions of Theorem (\ref{t1}) are satisfied if ${\mathbb E}l_{{\mathcal L}^{(1)}}<0$.
$l_{{\mathcal L}^{(1)}}$ can be estimated by
\[
-\nu \frac{1-\delta}{a_0\eps_{\cal L}^2}
+
\frac{2\nu b_0}{\nu\min\{h_1, h_2\}}(h_1\|\Delta\bar\psi_1\|_0^2+h_2\|\Delta\bar\psi_2\|_0^2)+2b_0\|\eta\|_0^2.
\]
On account of Lemma \ref{l0}
\[
\sup_{x\in B(\omega)}\frac{1}{t}\nu\int_0^t(h_1\|\Delta\bar\psi_1(\tau,\omega,x)\|_0^2+
h_1\|\Delta\bar\psi_2(\tau,\omega,x)\|_0^2)d\tau
\]
has the bound
\[
\frac{1}{t}\left(R(\omega)
+\int_0^t(d_0\|\eta(\theta_\tau\omega)\|_0^2
R(\theta_\tau\omega)+m(\theta_\tau\omega))d\tau\right)
\]
since $B$ is forward invariant.
Note that the expectation of $R/t$ can be made arbitrarily small if $t$ is large.
Calculation the expectation of the last expression and choosing $t$ sufficiently large
yields the conclusion.
\end{proof}

\begin{remark}\label{re-1} {\rm
i) It is easily seen that there exists a set $\hat{\mathcal{L}}^{h_1,h_2}$ of
linearly independent linear bounded functionals $\{\hat l_j\}_1^N$
which are determining with respect to $q$. Indeed, $q$ and $\psi$ are
connected by a linear homeomorphism  $\Lambda$ from ${\bf  L}_{per}^2$ to
${\bf H}_{per}^2$ defined by (\ref{Eqn3}) such that we can set
$\hat l^{h_1,h_2}_j=l^{h_1,h_2}_j\circ\Lambda$
where $\hat{\mathcal{L}}=\{l_j\}_1^N$ defines the set of determining functionals introduced
in Theorem \ref{main-1}.\\
ii) If $h_1=h_2$ then the set ${\cal L}$ itself is determining in the sense of Definition
\ref{dfp}.\\
iii)
Assumption (\ref{assump}) holds, if
\begin{equation}\label{assump-2}
\frac{4d_0 a_0 {\rm tr}_0Q}{\lambda_1^2\nu^2 (k+1)}<1,\quad
k+1>\frac{4}{a_0},
\end{equation}
for example. In this case it is easy to see that
$\left({\mathbb E} R_0^2\right)^{1/2}\le
4d_4^2 \left(\frac{3a_0}{2\nu\lambda_1}\right)^{2}$ (we choose $a=4/3$
in Lemma~\ref{l5}). Therefore using (\ref{constanten}) we obtain
\begin{equation*}
\Sigma\le C_4^{1/2}\frac{{\rm tr}_0Q}{\nu (k+1)}
\cdot\left(
4 d_0  d_4^2 \left(\frac{3a_0}{2\nu\lambda_1}\right)^{2}
+\nu\min\{h_1,h_2\} \right) +d_4.
\end{equation*}
Using (\ref{assump-2}) again we have the estimate
\begin{equation*}
d_4\le  \frac{{\rm tr}_0Q}{k+1}
\cdot\left( C_8^{1/4} \frac{d_1}
{4d_0 a_0} +C_4^{1/4} \frac{d_2}
{\lambda_1\nu}\right) +d_3.
\end{equation*}
Therefore under conditions (\ref{assump-2}) in the deterministic limit
${\rm tr}_0Q\to 0$ we obtain estimate

\begin{equation}\label{si-5}
\eps_{\cal L}<
\frac{\nu^2\lambda_1 \sqrt{\min\{h_1, h_2\}}}{6\|f\|_{-1} c_0\sqrt{h_1}}\cdot
\left\{\left(\lambda_1+2\max\{F_1, F_2\}\right)
\left(1+\frac{F_1 F_2}{\lambda^2_1}\right)\right\}^{-1/2}.
\end{equation}
If this estimate holds, then functionals $\{l^{h_1,h_2}_j\}$
are determining for
the {\em deterministic} two-layer model. Moreover under condition
(\ref{si-5}) these functionals are also determining for our
stochastic  two-layer model (\ref{q1}) and (\ref{pv1}) provided
the noise parameter ${\rm tr}_0Q$ is small enough.
}
\end{remark}

At the end of this section, we will discuss appropriate  values
for fluid parameters. We will estimate  the sufficient
condition  in Theorem \ref{main-1}, namely  the inequality (\ref{si-5}),
in terms of fluid parameters.

We further show that the long-time dynamics of the two-layer geophsysical
fluid system is
determined by the long time dynamics   of   the top layer  alone ,
when the fluid parameters and the Wiener process satisfy
certain conditons.  In this case the long time dynamics  will also
be  determined by only finitely many
functionals.\\
Let ${\cal L} =  \{ l_j : j= 1,...,N\}$ be a finite set of
linearly independent continuous functionals on the space $H^2_{per}$
of the top stream functions.
We assume that
\begin{equation}
 \| \nabla \psi_1 \|_0 \leq
\eps_{\cal L} \cdot  \|\D \psi_1\|_0
+  C_{\cal L} \cdot \max_{j= 1,\ldots, N} \vert l_j (\psi_1) \vert,
    \quad \psi_1 \in H^2_{per} ,
\label{condition}
\end{equation}
where $C_{\cal L}>0$ is a constant depending on ${\cal L}$.
As above we have
\[
h_1\|\D\psi_1\|_0^2 \ge \frac{1-\delta}{\eps_{{\cal L}}^2}h_1\|\nabla \psi_1\|_0^2
-C_{\delta,{\cal L}}h_1\max_{j=1,\cdots,N}|l_j(\psi_1)|^2.
\]
for any $0<\delta<1$ with appropriate positive constant
$C_{\delta,{\cal L}}$.

Thus under the condition (\ref{condition}) and using
the inequality
$\|\nabla \psi_2\|_0^2  \leq \lambda_1\|\D \psi_2\|_0^2$,
the estimate (\ref{estimate1a}) for the
solution difference $\hat{q}-\bar{q}$   becomes
\begin{equation*}
\begin{split}
\frac{d}{dt}\|q\|_{*}^2
&+ \nu h_1 \frac{1-\delta}{\eps_{{\cal L}}^2}\|\nabla \psi_1\|_0^2
+ (\nu \lambda_1 +2r) h_2  \|\nabla \psi_2\|_0^2
  \nonumber \\
&\leq
b_0\max\left\{\|\bar{\psi}_1-\xi_1\|_2^2\;,
\|\bar{\psi}_2-\xi_2\|_2^2\right\}\|q\|_{*}^2  +
C_{\delta,{\cal L}}\max_{j=1,\cdots,N}|l_j(h_1\psi_1)|^2.
\end{split}
\end{equation*}
Thus we obtain
\begin{equation*}
 \frac{d}{dt}\|q\|_{*}^2 \le l^{(2)}_{\cal L}(\bar{\psi}, \om)
\|q\|_{*}^2 + C_{\delta,{\cal L}}\max_{j=1,\cdots,N}|l_j(h_1\psi_1)|^2,
\end{equation*}
where
\begin{equation}
\begin{split}
 l^{(2)}_{\cal L}(\bar{\psi}, \om) =&
-\min\left\{ \nu \frac{1-\delta}{a_0\eps_{\cal L}^2}\; ,
\frac{\nu \lambda_1 +2r}{a_0}
\right\} \nonumber \\
&+
b_0\max\left\{\|\bar{\psi}_1-\xi_1\|_2^2\;,
\|\bar{\psi}_2-\xi_2\|_2^2\right\}.
\end{split}
\end{equation}
Again applying Theorem \ref{t1} we obtain the main result.

\begin{theorem}  \label{main}
Let    ${\cal L} =  \{ l_j : j= 1,...,N\}$ be a finite set of
linearly independent continuous functionals on the space $H^2_{per}$.
We assume that this set satisfies the following condition involving
only with the top fluid layer dynamical variable, i.e.,  the stream function
$\psi_1(x, y, t)$
\begin{equation}\label{com-def}
 \| \nabla \psi_1 \|_0 \leq
\eps_{\cal L} \cdot  \|\D \psi_1\|_0
+  C_{\cal L} \cdot \max_{j= 1,\ldots, N} \vert l_j (\psi_1) \vert,
    \quad \psi_1 \in H^2_{per}.
 \end{equation}
If
\begin{equation}\label{si-6}
\Sigma
<
\min\left\{ \frac{\nu}{\eps_{\cal L}^2}\; ,
\left(\nu \lambda_1 +2r\right)
\right\}\frac{\nu\min\{h_1,h_2\}}{2a_0b_0}\;,
\end{equation}
where $\Sigma$ is given by (\ref{si-2}),
then ${\cal L}$ is a set of asymptotically determining functionals in
probability for the  stochastically forced
two-layer quasigeostrophic fluid system
(\ref{Eqn1}), (\ref{Eqn3}).
\end{theorem}
The proof of this theorem is   the same as the proof of Theorem
\ref{main-1}.\\
\medskip\par

Relation (\ref{si-6}) holds if
\begin{equation}\label{si-7}
 \frac{\nu}{\eps_{\cal L}^2}\ge
\nu \lambda_1 +2r
\end{equation}
and
\begin{equation}\label{si-8}
\Sigma<
\left(\nu \lambda_1 +2r
\right)\frac{\nu\min\{h_1,h_2\}}{2b_0a_0}\;.
\end{equation}
The parameter $\Sigma$ depends on $r$ via $d_2$. Therefore it is not
clear whether (\ref{si-8}) holds for some $r$. However as in Remark~\ref{re-1}
 in the deterministic limit
${\rm tr}_0Q\to 0$ the estimate (\ref{si-8}) turns into the relation
\begin{equation}\label{si-9}
\nu \lambda_1 +2r>
\frac{36c^2_0\| f\|_{-1}^2h_1}{\nu^3\lambda_1 \min\{h_1,h_2\}}
\cdot
\left(1+\frac{2}{\lambda_1}\max\{F_1, F_2\}\right)
\left(1+\frac{F_1 F_2}{\lambda^2_1}\right)\;.
\end{equation}
This observation leads to the following assertion.
\begin{coro}  \label{main-c} Assume that (\ref{si-9}) holds.
Let    ${\cal L} =  \{ l_j : j= 1,...,N\}$ be a finite set of
linearly independent continuous functionals on the space $H^2_{per}$.
We assume that this set satisfies (\ref{com-def}) with the parameter
$\eps_{\cal L}$ satisfying (\ref{si-7}).
Then there exists $\gamma>0$ such that
${\cal L}$ is a set of asymptotically determining functionals in
probability for the  stochastically forced
two-layer quasigeostrophic fluid system (\ref{Eqn1}), (\ref{Eqn3})
provided ${\rm tr}_0Q\le \gamma$.

Consequencely, the asymptotic  probabilistic dynamics of the
stochastically forced two-layer quasigeostrophic fluid system
is determined only by the top fluid layer.
\end{coro}

\begin{remark}{\rm
Note that the main task is to prove the existence of determining functionals
of the main equation (\ref{q1}). It is easily seen by the structure
the transformations (\ref{zz}) a set ${\mathcal L}$ is determining
in probability for (\ref{Eqn1})
if and only if the same set ${\mathcal L}$ is determining
in probability for (\ref{q1}); see \cite{CheDuaSchm00}.}
\end{remark}
As Theorems~\ref{main-1} and \ref{main} show the problem
of describing of finite families of determining functionals is reduced
to the study of sets of functionals for which the estimate (\ref{2dieze})
(or (\ref{condition})) holds with $\eps_{\cal L}$ small enough.
It is also important to calculate the best possible value for
 $\eps_{\cal L}$ for the given family of functional. For this parameter
there is the estimate from below depending only  on a number of functionals
\cite{C3} and this estimate coinside with $\eps_{\cal L}$ when functionals
are modes (see Example~1 below). We refer to \cite{C3,C4}
for further discussion concerning an optimal choice of families
of functionals with properties like
(\ref{2dieze}) or (\ref{condition}).
We also note that the deterministic counterparts of
Theorems~\ref{main-1} and \ref{main} were proved in \cite{BC} for other
boundary conditions.

\medskip
{\it Example 1 (Determining modes)}
Let $\{ e_i(x,y)\}_{i=1}^\infty$ be the basis of eigenfunctions of $-\D$
with the periodic boundary conditions in $O$ such that the corresponding
eigenvalues possesses the property
\begin{equation}\label{eig-val}
0<\lambda_1\le \lambda_2\le \ldots.
\end{equation}
We note that every eigenfunction has the form
$L^{-1}\exp\{\frac{2\pi}{L}(j_1x+j_2y)\}$. However we numerate them such that
(\ref{eig-val}) holds. Suppose ${\cal L}_N =  \{ l_j : j= 1,...,N\}$ is the
set of the functionals on $H^2_{per}$ of the form
\[
l_j(u)=\int_O u(x,y)e_j(x,y) dO,\quad j=1,\ldots, N.
\]
Then one can prove (see, e.g. \cite{C3}) that estimate (\ref{condition})
holds for ${\cal L}_N$ with $\eps_{{\cal L}_N}=\lambda_{N+1}^{-1/2}$ and this
value is the best possible among all families of functionals consisting
of $N$ elements. Since $\lambda_{N}\sim c_0 N L^{-2}$
for $N$ large enough, we have that
$\eps_{{\cal L}_N}\sim \tilde{c}_0  L/\sqrt{N}$ with some absolute constant
$\tilde{c}_0$. Therefore we can use Theorem~\ref{main} to estimate number
of determining modes. We also note that using the family ${\cal L}_N$
of the functionals on $H^2_{per}$ we can easily  construct
a family $\tilde{\cal L}_N$
of the functionals on ${\bf H}^2_{per}=H^2_{per}\times H^2_{per}$
such that (\ref{2dieze}) holds with $\eps_{{\cal L}_N}=\lambda_{N+1}^{-1/2}$.

{\it Example 2 (Determining nodes)}
Let us consider the nodes
\[
(x_i, y_j)=\frac{L}{\sqrt{N}}\cdot (i, j)\in O, \quad 1\le i,j\le\sqrt{N},
\]
and define functionals $l_{ij}$ on  $H^2_{per}$ as  $\delta$-functions
at $(x_i, y_j)$, i.e by the formulas $l_{ij}(u)=u(x_i, y_j)$.
Let ${\cal L}=\{ l_{ij} \}$.
One can prove (see, e.g. \cite{C3} or \cite{JT2})
that estimate (\ref{condition})
holds for ${\cal L}$ with
$\eps_{{\cal L}}= \bar{c}_0  L/\sqrt{N}$
where $\bar{c}_0=\frac12$ is an absolute constant \cite{JT2}.
Thus $\eps_{{\cal L}}$ has the same order for large $N$ as
 is the best possible value $\eps_{{\cal L}_N}$ for families of functionals
consisting of $N$ elements.\\

We also note that the estimates for the completeness defect $\eps_{\cal L}$
and the Ekman constant $r$ given  in Theorems \ref{main-1} and \ref{main}
are rather crude.
    In the two-layer quasigeostrophic model considered above,
    the parameter $\nu$ is the molecular
    viscosity. For ocean water, $\nu = 10^{-6} m^2s^{-1}$.
   However, fluid turbulence at small scales  can act as an extra
    dissipative mechanism, thus calling
    for the substitution of the molecular viscosity
    by a much larger eddy viscosity. For example, in the two-layer quasigeostrophic
    model
    flow simulation in \cite{Oz}, the eddy viscosity  (we still use the
same notation as the molecular viscosity) is  taken as
$\nu=50 m^2s^{-1}$.
The gravitational acceleration $g=9.81 ms^{-2}$.

    At mid-latitude ($45 \; \mbox{degree} N$),  $f_0=8 \times 10^{-5} s^{-1}$ and
$\beta = 2.3\times 10^{-11} m^{-1}s^{-1}$.

    For large-scale flows at mid-latitude,  such as    the Gulf
    Stream in the Atlantic ocean,  the horizontal spatial scale $L$
  is at the order of $1000 km$.
Moreover,   layer depth $h_1$ and $h_2$ are
at the order of $500m$ each  for large scale flows
such as the Gulf Stream at the mid-latitude in the Atlantic ocean.
For the eddy viscosity $50 m^2s^{-1}$
   in \cite{Oz},  the   Ekman constant $r$ is
 then
at the order of $10^{-5} s^{-1}$. The number $r$ is large when
the eddy viscosity  is taken to be large.

The ocean water mean density $\rho_0$ is about
$1025 kg m^{-3}$ or $1.025 g cm^{-3}$.
$\rho_1$ and $\rho_2$ are in the vicinity of
$\rho_0$. The density difference $\rho_2-\rho_1\approx 25 kg m^{-3}$.
But it is this small density difference
that in turn determines pressure differences
and thereby drive the ocean circulation \cite{Washington}.

The mean wind   forcing $f$, i.e., the deterministic part of the
curl of the wind stree on the top fluid layer is usually taken  as
a stationary, or being time-averaged and even also zonally averaged,
sinusoidal function.  For example \cite{Dymnikov},
\[
f= \frac{2\pi \tau_0}{\rho_0 h_1L} \sin\frac{2\pi y}{L},
\]
where the wind tension $\tau_0$ is of order 1 $dyne/cm^2$ or
of order 0.1 $N/m^2$.
With this mean wind forcing and physical parameters
specified above,
 the inequality (\ref{si-5}) turns into the estimate
$\eps_{\cal L} \leq 1.35\cdot 10^{-3} m$
and therefore in the case of Examples 1 or 2
for the number $N$ of functionals we obtain the estimate
$N\approx 10^{18}$. Thus Theorem~\ref{main-1}
should be only considered as an qualitative assertion about
finite-dimensionality
of the long-time behaviour of the
stochastically forced two-layer quasigeostrophic fluid system.
As for Theorem~\ref{main}, a similar calculations show that
 the condition (\ref{si-9}) can be  valid under some special
choice of parameters. Therefore  Theorem~\ref{main} only predicts the
 possibility of a situation when the bottom fluid layer is slaved
by the top layer.

\section{Summary}

We have considered  asymptotic  probabilistic dynamics of the
stochastically forced two-layer quasigeostrophic fluid system.
We first transformed this system into a coupled system of
random partial differential equations and then show that
the asymptotic  probabilistic dynamics of this system depends
only on the top fluid layer, provided that the Wiener process
and the fluid parameters satisfy a certain condition, i.e., the
inequality (\ref{si-6});  see Theorem \ref{main}.
In particular, this latter condition is satisfied when
the trace of the covariance operator of the Wiener process is controled by
a certain upper bound (see Corollary \ref{main-c}) and
  the Ekman constant $r$ is sufficiently large
(see  the inequality (\ref{si-9})).
Note that the generalized time derivative of the Wiener process
models the fluctuating part of the wind stress forcing on the top
fluid layer, and the Ekman constant $r$ measures the     rate for
vorticity decay due to the friction in the bottom Ekman layer.


\bigskip

{\bf Acknowledgement.}
A part of this work was done at the
Oberwolfach Mathematical Research Institute, Germany,  while J.
Duan and B. Schmalfu{\ss}  were Research in Pairs
Fellows, supported by {\em Volkswagen Stiftung}.
J. Duan would like to thank  Tamay Ozgokmen, University of Miami,
for helpful discussions.
This work was partly supported by the NSF Grant
DMS-9973204.



\begin{thebibliography}{99}

\bibitem{Arn98}
L.~Arnold,
\newblock {\em {R}andom {D}ynamical {S}ystems}.
\newblock Springer, Berlin, 1998.


\bibitem{Arn00}
L.~Arnold,
\newblock Hasselmann's programm visited: {T}he analysis of stochasticity in
  deterministic climate models.
\newblock Report 450, Universit{\"a}t Bremen, Institut f{\"u}r Dynamische
  Systeme, 2000.

\bibitem{Benilov} E. S. Benilov,  On the stability of large-amplitude geostrophic flows in a
two-layer fluid:  the case of "strong" beta-effect.
{\em J. Fluid Mech.} {\bf 284} (1995), 137--158.


\bibitem{Berloff} P. Berloff and S. P. Meacham,
On the stability of equivalent-barotropic and baroclinic models
of the wind-driven circulation, preprint, 1999.

\bibitem{Berloff2} P. Berloff and J. C. McWilliams,
Large-scale, low-frequency variability in wind-driven ocean
gyres,  {\em J. Phys. Oceanogr.} {\bf 29} (1999),
 1925--1949.

\bibitem{Bernier94} C. Bernier,
{\it Existence of attractor for the quasi-geostrophic
approximation of the Navier-Stokes equations and estimate of its dimension.} Advances in Mathematical Sciences and Applications, {\bf 4} (2), pp 465-489 (1994).


\bibitem{BC}
Ch. Bernier-Kazantsev and I.D. Chueshov,
The finiteness of determining degrees of freedom for the
 quasi-geostrophic multi-layer ocean model,
{\em Nonlinear Analysis} {\bf 42} (2000), 1499-1512.




\bibitem{BF}
L.~Berselli and F.~Flandoli.
\newblock Remarks on determining projections for stochastic dissipative
  equations.
\newblock {\em Discrete and Continuous Dynamical Systems}, 5(8):197--214, 1999.




\bibitem{Brannan} J. Brannan, J. Duan and T. Wanner,
 Dissipative Quasigeostrophic Dynamics under Random Forcing ,
 {\em  J. Math. Anal. Appl.}  {\bf 228} (1998),   221--233.



\bibitem{C3}
I.~Chueshov,
\newblock Theory of functionals that uniquely determine asymptotic dynamics of
  infinite-dimensional dissipative systems.
\newblock {\em Uspekhi Mat. Nauk}, 53(4):77--124, 1998.
\newblock (in Russian). English translation in {\it Russian Mathematical
  Surveys}  53:731--776, 1998.


\bibitem{C4}
I.~Chueshov,
\newblock {\em Introduction to the Theory of Infinite-Dimensional Dissipative
  Systems}.
\newblock Acta, Kharkov, 1999.
\newblock (in Russian).


\bibitem{C5}
I.~Chueshov,
\newblock On determining functionals for stochastic {N}avier - {S}tokes
  equations.
\newblock {\em Stochastics and Stochastics Reports}, 68:45--64, 1999.


\bibitem{CheDuaSchm00}
I.~Chueshov, J. Duan and B.~Schmalfu{\ss}.
\newblock Determining functionals for random partial differential equations,
\newblock  To appear in
{\em Nonlinear Differential Equations and Applications}, 2001.

\bibitem{Constantin} P. Constantin and C. Foias,
{\em Navier-Stokes Equations}, Univ. of Chicago Press,
Chicago, 1988.


\bibitem{DelSole-Farrell} T. DelSole and B. F. Farrell,
A stochastically excited linear system as a model for quasigeostrophic
turbulence: Analytic results for one- and two-layer fluids,
{\em J. Atmos. Sci.} {\bf 52} (1995) 2531-2547.


\bibitem{Dymnikov} V.  Dymnikov and E. Kazantsev,
On the genetic ``memory" of chaotic attractor of the
barotropic ocean model.
In Proceedings of the third bilateral conference
``Predictability of atmospheric and oceanic circulations "
of the French-Russian A.M.Liapunov Institute in Computer Science
 and Applied Mathematics (INRIA - Moscow State University). Nancy,
April, 1997.Edition MSU, 1997, pp. 25-36.


\bibitem{FL}
F.~Flandoli and J.~A. Langa.
\newblock On determining modes for dissipative random dynamical systems.
\newblock {\em Stochastics and Stochastics Reports}, 66:1--25, 1999.


\bibitem{FMTT}
C. Foias, O. Manley, R.Temam and Y.M. Treve,
Asymptotic analysis of the Navier-Stokes equations,
{\it Physica D} {\bf 9} (1983), 157--188.


\bibitem{FP}
C.~Foias and G.~Prodi.
\newblock Sur le comportement global des solutions nonstationnaires des
  \'{e}quations de {N}avier-{S}tokes en dimension deux.
\newblock {\em Rend. Sem. Mat. Univ. Padova}, 39:1--34, 1967.


\bibitem{FTiti} C. Foias, E.S. Titi,  Determining nodes, finite difference schemes and inertial manifolds. {\em Nonlinearity}, {\bf 4}, pp 135-153 (1991).







\bibitem{Griffa} A. Griffa and S. Castellari,
Nonlinear general circulation of an ocean model driven by
wind with a stochastic component,
{\em J. Marine Research} {\bf 49} (1991), 53-73.

\bibitem{Has76}
K.~Hasselmann.
\newblock Stochastic climate models, {P}art {I}.
\newblock {\em Tellus}, 28:473--485, 1976.


\bibitem{Holloway} G. Holloway,
Ocean circulation:  Flow in probability under
statistical dynamical forcing,
in {\em Stochastic Models in Geosystems},
S. Molchanov and W. Woyczynski (eds.), Springer, 1996.

\bibitem{Huang} R. X. Huang and H. Stommel,
Cross sections of a two-layer inertial Gulf Stream,
{\em J. Phys. Oceanography} {\bf 20} (1990), 907-901.

\bibitem{JT2} D.A. Jones and  E.S. Titi,
Determining  finite volume elements for the 2D Navier-Stokes equations,
{\em Physica D} {\bf 60} (1992), 165-174.






\bibitem{La2}
O.~Ladyzhenskaya.
\newblock A dynamical system generated by the {N}avier--{S}tokes equations.
\newblock {\em Journal of Soviet Mathematics}, 3:458--479, 1975.



\bibitem{Muller} P. M\"uller, Stochastic forcing of quasi-geostrophic eddies,
in {\em Stochastic Modelling in Physical Oceanography},
R. J. Adler, P. M\"uller and B. Rozovskii (eds.), Birkh\"auser, 1996.

\bibitem{Oz} T. Ozgokmen,
Emergence of inertial gyres in a two-layer
quasigeostrophic ocean model,
{\em J. Phys. Oceanography} {\bf 28} (1998), 461-484.




\bibitem{Ped87}
J.~Pedlosky.
\newblock {\em Geophysical Fluid Dynamics}.
\newblock Springer Verlag, New-York, Berlin, 1987.


\bibitem{Ped96}  J.~Pedlosky,
 {\em Ocean Circulation Theory}.
Springer--Verlag,  Berlin, 1996.


\bibitem{Sal98}
R.~Salmon.
\newblock {\em Lectures on Geophysical Fluid Dynamics}.
\newblock Oxford Univ. Press, Oxford, 1998.

\bibitem{Salmon}
R.~Salmon.
Generalized two-layer models of ocean circulation,
{\em J. Marine Research} {\bf 52} (1994),  865-908.

\bibitem{Samelson} R. M. Samelson, Stochastically forced current
fluctuations in vertical shear and over topography,
{\em J. Geophys. Res.} {\bf 94} (1989) 8207-8215.
\bibitem{Washington} W. M. Washington and C. L. Parkinson,
{\em An Introduction to Three-Dimensional Climate Modeling},
Oxford Univ. Press, 1986.

\end{thebibliography}
\end{document}